\input amstex
\documentstyle{amsppt}
\input xypic
\input psfig
\magnification=1200
%\pagewidth{5.4in}
\pageheight{7.6in}
\expandafter\redefine\csname logo\string@\endcsname{}
\NoBlackBoxes
%\NoRunningHeads

\topmatter
\title The Fibered isomorphism conjecture for complex manifolds 
\endtitle
\author S. K. Roushon
\endauthor
\thanks February 18, 2005 \endthanks
\endtopmatter
\centerline {School of Mathematics}
\centerline {Tata Institute}
\centerline {Homi Bhabha Road}
\centerline {Mumbai 400 005, India.} 
\centerline {{\rm email address:-} roushon\@math.tifr.res.in}

\medskip

\noindent
{\bf Abstract.} In this paper we show that the Fibered Isomorphism
Conjecture of Farrell and Jones corresponding to the stable topological
pseudoisotopy functor is true for the fundamental groups of a class of
complex manifolds. A consequence of this result is that the Whitehead
group, reduced projective class groups and the negative $K$-groups of the
fundamental groups of these manifolds vanish whenever the fundamental
group is torsion free. We also prove the same results for a class of real
manifolds including a large class of $3$-manifolds which has a
finite sheeted cover fibering over the circle.

\medskip



\medskip
\noindent
{\it Abbreviated title:} Isomorphism conjecture for complex manifolds

\document
\baselineskip 14pt

\newpage
\head 
0. Introduction
\endhead

In this paper we consider proving the Farrell-Jones Fibered Isomorphism
Conjecture (FIC) corresponding to the stable topological pseudoisotopy
functor for the fundamental groups of complex manifolds. 
This conjecture is already proved for several classes of groups
including fundamental groups of closed nonpositively curved Riemannian
manifolds (\cite{6}), for cocompact discrete subgroup of virtually
connected Lie groups (\cite{6}) and for the class of virtually strongly
poly-free groups (\cite{9}). A rich class of complex manifolds
are smooth complex projective algebraic varieties. Some results are 
known on the structure of the fundamental groups of a large class of
complex surfaces. We make use of these informations for the proofs of
our results. The method of the proofs also generalizes to consider
some real manifolds and to prove the FIC for the fundamental groups. 

The main input to the proofs of Theorem 1.3 and 1.5 comes from
theorem 7.1 of Farrell and Linnell in \cite{7} where they proved that if
the Fibered Isomorphism Conjecture is true for all the groups in a
directed system of groups then it is true for the direct limit also. The
proof of our Main Lemma uses this result to show that the FIC is
true for a large class of mapping torus of the fundamental group of a 
closed orientable real $2$-manifold and also for a certain class of
mapping tori of infinitely generated free groups. The key idea to prove
the Main Lemma was that except for two closed surfaces the covering space
corresponding to the commutator subgroup of the fundamental group of all
other surfaces have one topological end. Though the spaces involved in the
theorems have finitely presented fundamental groups, during the proof we
encounter some infinitely generated groups and these are the places
where we use the Main Lemma crucially.

Throughout the paper whenever we encounter a $3$-manifold which fiber
over the circle with fiber a surface of genus $\geq 2$ we make the
assumption that the monodromy diffeomorphism
is special (see Definition in Section 1). In \cite{27} we have proved that
we can remove this assumption of being `special' provided the FIC is true
for $A$-groups (see \cite{27} for definition). Roughly speaking a 
torsion free $A$-group is a discrete group which is isomorphic to the
fundamental group of a complete nonpositively curved Riemannian manifold
whose metric is $A$-regular. In a recent paper (\cite{16}) L.E. Jones
proved that the assembly map in the statement of the FIC induces a
surjective homomorphism on the homotopy group level for any torsion free
$A$-group. He also stated some conjectured theorems ([\cite{16},
conjectured theorems 6.7 and 6.9]) which imply that the
`torsion free' assumption can be dropped from the statement of the result.
On the other hand the FIC states that the above homomorphism should be an
isomorphism.

In Section 1 using Lefschetz hyperplane section theorem we deduce that if 
the FIC is true for the fundamental group of smooth 
complex algebraic surfaces then it is true for the fundamental group
of any smooth complex algebraic variety. Also we state our results in this
section. Section 2 recalls the Farrell and Jones Fibered Isomorphism
Conjecture and states the known results we need. In Section 3 we make a
brief trip to classification of complex surfaces and their topological
properties. Sections 4 and 5 contain the proofs of the results stated in
Section 1 and we prove the Main Lemma in Section 6. In Section 7 we prove
the FIC for a class of virtually fibered $3$-manifolds. Section 8
contains examples of special diffeomorphisms.

\head
1. Reduction to surface case and statements of results
\endhead

Let $X$ be a smooth complex projective algebraic variety of dimension $n$.
By definition $X\subset {\Bbb {CP}}^m$ for some $m$ and is a 
complex submanifold of ${\Bbb {CP}}^m$. It is well-known that not all
complex manifold is a complex submanifold of some complex
projective space; otherwise it will become algebraic. In Section 3 we will
mention such examples in the case of surfaces. 

There is a natural collection of complex submanifolds of $X$ arising 
from taking intersection of $X$ with hyperplanes $H$ in ${\Bbb {CP}}^m$.
For a general hyperplane $H$ the intersection of $H$ with $X$ is a  
connected complex submanifold of $X$ ([\cite{2}, chapter I, corollary
20.3]). Let $H_0$ be such a hyperplane.
Then the two manifolds $X$ and $X\cap H_0$ shares similar 
homological and homotopical properties up to a certain degree. This is the
content of the Lefschetz hyperplane section theorem ([\cite{2}, chapter
I, theorem 20.4]). 

\proclaim{Lefschetz hyperplane section theorem} For $n\geq 2$ the
inclusion map $X\cap H_0\subset X$ induces the following isomorphisms.
$$H_i(X\cap H_0, {\Bbb Z})\to H_i(X, {\Bbb Z})$$ $$\pi_i(X\cap H_0, {\Bbb
Z})\to \pi_i(X, {\Bbb Z})$$

for $0\leq i\leq n-2$.\endproclaim

Let ${\Cal G}_n$ be the class of fundamental groups of all smooth complex 
projective algebraic varieties of dimension $n$. Then successively applying
Lefschetz hyperplane section theorem we get the following Lemma.

\proclaim{Lemma 1.1} $\cup_{n=2}^{\infty}{{\Cal G}_n}\subset 
{\Cal G}_2$.\endproclaim

Thus from Lemma 1.1 we see that if we want to prove the FIC for the
fundamental group of smooth projective algebraic varieties then it is
enough to consider the smooth projective algebraic surfaces. 

The rest of this section contains the statements of the results. Before
that we recall some standard definition from algebraic geometry. By a {\it
complex surface} we mean a closed complex $2$-manifold. By an
{\it algebraic} complex surface we mean it is a complex surface and 
is defined by finitely many homogeneous polynomials in $n+1$ variable in
the complex projective space ${\Bbb {CP}}^n$. By a {\it curve} we will 
mean a complex projective algebraic variety of dimension $1$. $\kappa(X)$
denotes the Kodaira dimension of $X$. For definition of $\kappa(X)$ see
[\cite{2}, p. 23] or [\cite{10}, definition 1.6]. When $X$ is a complex
surface $\kappa(X)\in \{-\infty,0,1,2\}$.

\proclaim{Theorem 1.2} Let $X$ be a complex surface of one of the 
following types.

\roster
\item $X$ is algebraic and $\kappa(X)=-\infty$

\item $X$ is a Hopf surface

\item $\kappa(X)=0$

\item $X$ is an Inoue surface 

\item $X$ is an elliptic surface

\endroster

Then the Fibered Isomorphism Conjecture is true for
$\pi_1(X)$.\endproclaim

The theorem below deals with some more complex manifolds. To state the
theorem we need to make some definition. At first recall that it follows
from a result of Hillman [\cite{12}, theorem 7] that if a complex surface
fibers over the circle then the fiber is a Seifert fibered space (see
Theorem 4.1) . Hence if $X$ is such a surface then $\pi_1(X)\simeq
\pi_1(S)\rtimes {\Bbb Z}$ where $S$ is a Seifert fibered space. Assume that
the monodromy diffeomorphism is a fiber preserving diffeomorphism of the
Seifert fibered space $S$. If $\pi_1(S)$ is infinite then there is an
infinite cyclic normal subgroup of $\pi_1(S)$ with quotient
$\pi_1^{orb}(B)$ where $B$ is the base orbifold of $S$ ([\cite{11},
chapter 12]). Again if
$\pi_1^{orb}(B)$ is infinite then one can find a finite index
characteristic subgroup $K$ of $\pi_1^{orb}(B)$ which is isomorphic to a
closed surface group (see Section 4). Note that in this situation the
monodromy diffeomorphism of the fiber bundle $X\to {\Bbb S}^1$ induces an
automorphism of $\pi_1^{orb}(B)$. Since $K$ is characteristic we have an
exact sequence $$1\to K\to \pi_1^{orb}(B)\rtimes {\Bbb Z}\to
(\pi_1^{orb}(B)/K)\rtimes {\Bbb Z}\to 1$$   

Let $l$ be an element of $(\pi_1^{orb}(B)/G)\rtimes {\Bbb Z}$ which
generates an infinite cyclic normal subgroup of finite index. Since $K$
is a closed surface group the action of $l$ (by conjugation by a lift of 
$l$) on $K$ is induced, up to conjugation, by a diffeomorphism $f_l$ of a
closed surface $F$ so that $\pi_1(F)$ is isomorphic to $K$.
Let us call $f_l$ a {\it base} diffeomorphism associated to the infinite
cyclic normal subgroup generated by $l$. 

\proclaim{Theorem 1.3} Let $X$ be a complex surface which is the total
space of a fiber bundle over the circle ${\Bbb S}^1$. Under the above
notations when $F$ is not the $2$-torus assume that there is a base
diffeomorphism $f_l$ which is special (see Definition below). Then the
FIC is true for $\pi_1(X)$.\endproclaim  

In Section 4 we will also prove that the FIC is true for a class of
complex
surfaces of Kodaira dimension $2$. A large class of examples of surfaces
of Kodaira dimension $2$ are ramified $2$-sheeted covering of ${\Bbb
{CP}}^2$ ramified along a curve. We will prove the FIC for a class of such
surfaces and will give an example to show that, given the methods
available, this is the best possible result we can prove.  

Complex surfaces with Kodaira dimension $1$ are elliptic surfaces (Theorem
3.3.1). There is a notion of elliptic surfaces in the smooth category
called $C^{\infty}$-elliptic surface. These are smooth real $4$-manifold
which are locally modelled on complex elliptic surfaces. We recall the
definition of $C^{\infty}$-elliptic surface in Section 3. The
fundamental groups of these $4$-manifolds have close properties with
that of complex elliptic surfaces. We prove in Corollary 1.4 that the
FIC is true for a class of these manifolds also. 

\proclaim{Corollary 1.4} The FIC is true for the fundamental group of a
$C^{\infty}$-elliptic surface $X$ if $X$ has no singular fiber and with
cyclic monodromy.\endproclaim

There is another natural collection of
smooth $4$-manifolds which are fiber bundles over real $2$-manifolds with
real $2$-manifolds as fiber. In Theorem 1.5 below we prove that the
FIC is true for the fundamental group of a large class of manifolds from
this collection. A large class of complex surfaces belong
to this collection where both the fiber and base are $2$-manifolds of
genus $\geq 2$. In this particular case the fiber bundle projection is
called {\it Kodaira fibration} (see [\cite{2}, chapter V, section 14] 
for details). To state our next results we need the following definition.

\proclaim{Definition} {\rm Let $F$ be a closed orientable surface and $f$
is an orientation preserving diffeomorphism of $F$. Let $\tilde F$ be the
covering of $F$ corresponding to the commutator subgroup of $\pi_1(F)$ and
let $\tilde f$ be the lift of $f$ to $\tilde F\to \tilde F$. Let 
$p:M_{\tilde f}\to M_f$ be the covering projection from the mapping torus
of $\tilde f$ to that of $f$. 

We say $f$ is a {\it special} diffeomorphism if one of the following
holds.

\roster
\item the mapping torus of $f$ supports a nonpositively curved Riemannian
metric.
\item some power of $f$ is isotopic to identity.
\item the fundamental group of any component of $p^{-1}(S)$ is not free,
where $S$ varies over all Seifert fibered pieces in the Jaco-Shalen and
Johannson (JSJ) decomposition of the mapping torus $M(f)$.\endroster
\endproclaim

We now recall JSJ-decomposition of a $3$-manifold briefly. A $3$-manifold
is called {\it irreducible} if any embedded $2$-sphere in the manifold
bounds an embedded $3$-disc. A compact orientable irreducible $3$-manifold
$M$ is called {\it Haken} if there is a compact orientable surface $S$
embedded in $M$ such that $\pi_1(S)$ is infinite and the inclusion map
$S\to M$ induces an injective homomorphism $\pi_1(S)\to \pi_1(M)$. The
JSJ-decomposition states that any Haken $3$-manifold admits a
decomposition along finitely many mutually nonparallel tori embedded in
$M$ so that the decomposed pieces are either Seifert fibered or simple
(see \cite{14}, \cite{15}). Thurston proved that these simple pieces
admits complete hyperbolic metric in the interior. 

For examples of special diffeomorphisms, recall that the mapping
tori of pseudo-Anosov diffeomorphisms are hyperbolic (\cite{22},
\cite{28}). Also a large class of examples of special diffeomorphisms
satisfying condition $(3)$ above is given in Section 8.  

\proclaim{Theorem 1.5} Let $X$ be an orientable real $4$-dimensional
manifold which is the total space of a fiber bundle over an orientable
real $2$-dimensional manifold with fiber an orientable real
$2$-dimensional manifold. Assume all the monodromy diffeomorphisms of
the fiber are special. Then the FIC is true for $\pi_1(X)$.\endproclaim

In fact we prove that the FIC is true for any extension of a (real)
surface
group by a (real) surface group under a similar hypothesis.

From the proof of the above theorem the following more general corollary
is easily deduced.

\proclaim{Corollary 1.6} Let $M^{n+2}\to N^n$ be a fiber bundle projection
of real manifolds with fiber $F$. Let $\pi_1(N)$ be torsion free,
$\pi_2(N)=1$ and the 
FIC is
true for $\pi_1(N)$. If $M$ is not compact then assume in addition that
the fiber of the fiber bundle projection $M^{n+2}\to N^n$ is either of
finite topological type or $F$ has infinitely generated fundamental
group and $F$ is the covering of a compact surface $F'$ corresponding to
the commutator subgroup of $\pi_1(F')$ and any monodromy diffeomorphism of
the fiber $F$ is a lift of a special diffeomorphism of $F'$. Also 
assume the monodromy diffeomorphisms of the fiber $F$ are
special when $F$ is closed and of genus $\geq 2$. Then the FIC is true for
$\pi_1(M)$.\endproclaim

The main ingredient behind the proof of the above results is the following
proposition.

\proclaim{Proposition 1.7} Let $N$ be a closed orientable $3$-manifold and
there is a finite sheeted cover $M$ of $N$ which fibers over the circle. 
Assume that the monodromy diffeomorphism of the fiber bundle $M\to
{\Bbb S}^1$ is special when the fiber is a surface of genus $\geq 2$. Then
the FIC is true for $\pi_1(N)$.\endproclaim

An important corollary to the above results is the following.

\proclaim{Corollary 1.8} Let $G$ be a torsion free subgroup of
$\pi_1(X)$ where $X$ is a space appearing in the Theorems 1.2, 1.3 and
1.5, Corollaries 1.4 and 1.6 and Proposition 1.7. Then $Wh(G)={\tilde
K}_0({\Bbb Z}G)=K_{-i}({\Bbb Z}G)=0$ for all $i\geq 1$.\endproclaim

We end this section with the following remark.

\proclaim{Remark 1.9} {\rm We have already mentioned in the introduction
that throughout this section the assumption `special' for the monodromy
diffeomorphism of a fibered $3$-manifold can be dropped provided 
the FIC is true for $A$-groups. Here we remark that in fact we do not need
to assume this strong result. In \cite{27} we show that we only need to
assume that the FIC is true for the fundamental groups of compact
irreducible $3$-manifolds whose boundary components are all incompressible
and are surfaces of genus $\geq 2$ (we called these groups as $B$-groups).  
Also we showed that a $B$-group is an $A$-group. We also proved (in
\cite{27}) that the FIC is true for a large class of 
$B$-groups.}\endproclaim

\head
2. Farrell-Jones fibered isomorphism conjecture 
\endhead

In this section we recall the Fibered Isomorphism Conjecture of 
Farrell and Jones made in \cite{6}. 

Let $\Cal S$ denotes one of the three functors from the category of
topological spaces to the category of spectra: (a) the stable topological
pseudoisotopy functor ${\Cal P}()$; (b) the algebraic $K$-theory functor 
${\Cal K}()$; (c) and the $L$-theory functor ${\Cal L}^{-\infty}()$. 

Let $\Cal M$ be the category of continuous surjective maps. The
objects of $\Cal M$ are continuous surjective maps $p:E\to B$ between
topological spaces $E$ and $B$. And a morphism between two maps $p:E_1\to
B_1$ and $q:E_2\to B_2$ is a pair of continuous maps $f:E_1\to E_2$,
$g:B_1\to B_2$ such that the following diagram commutes. 

$$\diagram
E_1 \rto^f \dto^p & E_2 \dto^q\\
B_1 \rto^g &B_2\enddiagram$$

There is a functor defined by Quinn \cite{25} from $\Cal M$ to the
category
of $\Omega$-spectra which associates to the map $p:E\to B$ the spectrum
${\Bbb H}(B, {\Cal S}(p))$ with the property that ${\Bbb H}(B, {\Cal
S}(p))={\Cal S}(E)$ when $B$ is a single point. For an explanation of
${\Bbb H}(B, {\Cal S}(p))$ see [\cite{6}, section 1.4]. Also the map
${\Bbb H}(B, {\Cal S}(p))\to {\Cal S}(E)$ induced by the morphism:
id$:E\to E$; $B\to *$ in the category $\Cal M$ is called the Quinn
assembly map. 

Let $\Gamma$ be a discrete group and $\Cal E$ be a $\Gamma$ space which
is universal for the class of all virtually cyclic subgroups of $\Gamma$ 
and denote ${\Cal E}/\Gamma$ by $\Cal B$. For definition of universal 
space see [\cite{6}, appendix]. Let $X$ be a space on which $\Gamma$ acts
freely and properly discontinuously and $p:X\times_{\Gamma} {\Cal E}\to
{\Cal E}/{\Gamma}={\Cal B}$ be the map induced by the projection onto the 
second factor of $X\times {\Cal E}$. 

The Fibered Isomorphism Conjecture states that the map $${\Bbb
H}({\Cal B}, {\Cal S}(p))\to {\Cal S}(X\times_{\Gamma} {\Cal E})={\Cal
S}(X/\Gamma)$$ is a (weak) equivalence of spectra. The equality is induced
from the map $X\times_{\Gamma}{\Cal E}\to X/\Gamma$ and using 
the fact that $\Cal S$ is homotopy invariant.  

Let $Y$ be a connected $CW$-complex and $\Gamma=\pi_1(Y)$. Let $X$ be the
universal cover $\tilde Y$ of $Y$ and the action of $\Gamma$ on $X$ is
the action by group of covering transformation. If we take an aspherical
$CW$-complex $Y'$ with $\Gamma=\pi_1(Y')$ and $X$ is the universal cover
$\tilde Y'$ of $Y'$ then by [\cite{6}, corollary 2.2.1] if the FIC is
true for the space $\tilde Y'$ then it is true for $\tilde Y$ also. Thus
whenever we say that the FIC is true for a discrete group $\Gamma$ or for
the
fundamental group $\pi_1(X)$ of a space $X$ we would mean it is true for
the Eilenberg-MacLane space $K(\Gamma, 1)$ or $K(\pi_1(X), 1)$ and the
functor ${\Cal S}()$.  

Throughout this paper we consider only the stable topological
pseudoisotopy functor; that is the case when ${\Cal S}()={\Cal P}()$. 
And by the FIC we mean the FIC for ${\Cal P}()$.

Now we recall few known results we need about the FIC. 

\proclaim{Lemma A} ([\cite{6}, theorem A.8]) If the FIC is true for a
discrete group $\Gamma$ then it is true for any subgroup of
$\Gamma$.\endproclaim 

\proclaim{Lemma B} Let $\Gamma$ be an extension of  
the fundamental group $\pi_1(M)$ of a closed nonpositively curved 
Riemannian manifold or a compact orientable surface (may be with nonempty 
boundary) $M$ by a finite group $G$ then the FIC is true for
$\Gamma$. Moreover the FIC is true for the wreath product
$\Gamma\wr G$.\endproclaim

\proclaim{Lemma C} ([\cite{6}, proposition 2.2]) Let $f:G\to H$ be a
surjective homomorphism. Assume
that the FIC is true for $H$ and for $f^{-1}(C)$ for all virtually cyclic
subgroup $C$ of $H$ (including $C=1$). Then the FIC is true for
$G$.\endproclaim

We will use Lemma A, Lemma C and [\cite{9}, algebraic lemma] throughout
the paper, sometimes even without referring to it.

\proclaim{Lemma D} Let $S$ be a closed $2$-dimensional orbifold. Then
the FIC
is true for $\pi_1^{orb}(S)$.\endproclaim

Before we give the proofs of Lemma B and D we recall some group
theoretic definition. Let $G$ and $H$ be two groups. Assume $G$ is finite.
Then $H\wr G$ denotes the wreath product with respect to the regular
action of $G$ on $G$. Recall that actually $H\wr G\simeq H^G\rtimes G$
where $H^G$ is product of $|G|$ copies of $H$ indexed by the elements of
$G$ and the action of $G$ on the product is induced by the regular action
of $G$ on $G$.

\demo{Proof of Lemma B and D} Let us start with the hypothesis of
Lemma D. If $\pi_1^{orb}(S)$ is finite then there is nothing to prove
because the FIC is true for any finite group. 
So assume it is infinite. If the orbifold fundamental group
is infinite then there is a finite index subgroup $H$ of $\pi_1^{orb}(S)$
such that $H$ is isomorphic to the fundamental group of a closed surface.
Taking intersection of all conjugates of $H$ in $\pi_1^{orb}(S)$ we get a 
finite index normal subgroup $H_1 < H$ of $\pi_1^{orb}(S)$. Clearly
$H_1$ is again the fundamental group of a closed surface, say $\tilde
S$. Thus we have an exact sequence $$1\to \pi_1(\tilde S)\to
\pi_1^{orb}(S)\to G\to 1.$$ Here $G$ is a finite group. 

Let $\Gamma=\pi_1^{orb}(S)$ and $M=\tilde S$. We have reached the
hypothesis of Lemma B in the case when $M$ is closed.
 
By [\cite{9}, algebraic lemma] or [\cite{4}, theorem 2.6A] we have an
embedding of $\Gamma$ in the wreath product $\pi_1(M)\wr G$, where the
wreath product is taken using the regular action of the group $G$ on $G$. 
Let $U=M\times \cdots \times M$ be the $|G|$-fold product of $M$.
Then $U$ is a closed nonpositively curved Riemannian manifold. By
[\cite{9}, fact 3.1] it follows that the FIC is true for
$\pi_1(U)\rtimes G\simeq (\pi_1(M))^G\rtimes G\simeq \pi_1(M)\wr G$.
Lemma A now proves that the FIC is true for $\Gamma$.

If $M$ is a compact surface with nonempty boundary then $\pi_1(M) <
\pi_1(N)$ where $N$ is a closed nonpositively curved surface. Hence 
$\Gamma < \pi_1(M)\wr G < \pi_1(N)\wr G$. Using Lemma A we
complete the proof.\qed\enddemo  

Finally we recall the following important case when the FIC is true.

\proclaim{Theorem E} ([\cite{6}, theorem 2.1]) Let $\Gamma$ be a
cocompact discrete subgroup of a virtually connected Lie group. Then
the FIC is true for $\Gamma$.\endproclaim

We will also use the following consequence of [\cite{6}, proposition 2.4]
frequently. Recall that a poly-$\Bbb Z$ group is a group having a normal
series whose normal quotients are infinite cyclic. And a group is
virtually poly-$\Bbb Z$ if it has a finite index subgroup which is
poly-$\Bbb Z$. 

\proclaim{Theorem F} (\cite{6}) The FIC is true for any virtually
poly-$\Bbb Z$ group.\endproclaim 

\head
3. Classification of complex surfaces and related topological results 
\endhead

In this section we recall the classification of complex surfaces and their 
topological properties which we need in the next section for the proof of 
Theorems 1.2 and Corollary 1.4. The references for this material we follow
are \cite{10} and \cite{2}. 

By a {\it complex surface} we mean a compact complex manifold of
complex dimension $2$. There are two classes of complex manifolds; the
algebraic and the nonalgebraic ones. Throughout the paper by {\it
algebraic surface} we will mean (unless otherwise stated) a complex
surface which is a complex submanifold of some complex 
projective space ${\Bbb {CP}}^n$ and also this condition is equivalent to 
saying that it can be obtained as a set of zeros of finitely
many homogeneous polynomial in some complex projective space. 

Let $X$ be a complex surface. For a point $x\in X$ the {\it blow up}
of $X$ at $x$ is the surface $X\# {\overline {\Bbb {CP}}^2}$ where the
connected sum is taken around a ball at $x$. Here ${\overline {\Bbb
{CP}}^2}$ is ${\Bbb {CP}}^2$ with the opposite orientation. By blowing up
at a point $x$ we introduce a curve $C$ which is isomorphic to ${\Bbb
{CP}}^1$ in the space replacing the point $x$. $C$ has self-intersection
number $-1$. This is also called a $-1$ curve. Conversely if there is a
$-1$ curve $C$ in a complex surface $X'$ then there is another complex
surface $X$ and a point on it so that $X'$ is obtained from $X$ by blowing
up $X$ at the point $x$. This is called {\it blowing down} the curve $C$.
In this discussion our main interest is that by blowing down
or blowing up in a complex surface we do not change the fundamental group.    

For the definition of {\it Kodaira dimension} $\kappa(X)$ of a complex
surface $X$ we refer the reader to [\cite{10}, definition 1.6].
$\kappa(X)$ can assume only $4$ values; $-\infty, 0, 1$ and $2$. There is
a classification of complex surface in terms of this dimension and we
recall the known topological properties of surfaces of different Kodaira
dimension. 

We recall these results without proof and give the references where the
proof can be found.

\subhead 
3.1 The case $\kappa(X)=-\infty$
\endsubhead

A complex surface $X$ is called {\it ruled} ([\cite{10}, definition 1.9])
if there is a holomorphic map 
$f:X'\to S$, where $S$ is a complex $1$-manifold and all the fibers of 
$f$ are isomorphic to ${\Bbb {CP}}^1$ and $X$ is a blow up of $X'$. 
The map $f:X'\to S$ is called a
ruling of $X'$. It is well known that given a ruling $f:X\to S$ on a
complex surface there is a rank $2$ complex vector bundle on $S$ whose 
associated ${\Bbb {CP}}^1$ bundle is isomorphic to $X$ (see [\cite{10},
chapter I, section 1.2.1]). Thus we see that $f$
is in fact a fiber bundle projection. Hence for any ruled surface
$X$, $\pi_1(X)$ is isomorphic to $\pi_1(S)$ for some curve $S$. 

\proclaim{Proposition 3.1.1} ([\cite{10}, chapter I, theorem 1.10]) Let
$X$
be a complex algebraic surface with $\kappa(X)=-\infty$ then $X$ is
either ${\Bbb {CP}}^2$ or is ruled.\endproclaim 

There are nonalgebraic surfaces with $\kappa(X)=-\infty$. All of them has
first Betti number equal to one and with infinite fundamental group. A
classification of these surfaces are not yet known. A class of examples of
such surfaces are {\it Hopf surfaces} (see [\cite{10}, definition 7.13]).
These are complex surfaces which
has universal cover biholomorphic to ${\Bbb C}^2-(0,0)$. Also it is known
([\cite{10}, chapter I, proposition 7.15]) that any Hopf surface has a
finite sheeted cover homotopy equivalent to ${\Bbb S}^1\times {\Bbb S}^3$.
Hence we have the following Proposition.

\proclaim{Proposition 3.1.2} Any Hopf surface has virtually infinite
cyclic fundamental group.\endproclaim

Another class of example of nonalgebraic surfaces with
$\kappa(X)=-\infty$ are Inoue surfaces ([\cite{2}, chapter V, section
19]). These surfaces are by construction fiber bundle over ${\Bbb S}^1$
with fiber a $3$-manifold which is a principal ${\Bbb S}^1$-bundle over
the $2$-torus ${\Bbb S}^1\times {\Bbb S}^1$.

\proclaim{Proposition 3.1.3} \cite{23} The fundamental group of an Inoue
surface is isomorphic to $G\rtimes {\Bbb Z}$ where $G$ is the fundamental
group of a principal ${\Bbb S}^1$ bundle over the $2$-torus.\endproclaim

\proclaim{Corollary 3.1.4} The fundamental group of an Inoue surface is
poly-$\Bbb Z$.\endproclaim

\subhead
3.2 The case $\kappa(X)=0$
\endsubhead

The following two propositions give the required topological properties
we need of surfaces with $\kappa(X)=0$.

\proclaim{Proposition 3.2.1} ([\cite{10}, chapter I, theorem 2.6]) Let $X$
be an algebraic surface and $\kappa(X)=0$. Then $\pi_1(X)\simeq {\Bbb
Z}^4, \{1\}, {\Bbb Z}/{2{\Bbb Z}}$, or $\pi_1(X_{\text min})$. Here
$X_{\text min}$ is a hyperelliptic surface. In the last case we have an
exact sequence $$1\to \Lambda\to \pi_1(X)\to G\to 1$$ where $G$ is a
finite group and $\Lambda$ is a lattice in ${\Bbb C}^2$.\endproclaim

By definition a hyperelliptic surface is a finite quotient of a complex
torus ([\cite{10}, chapter I, section 1.1.4])

\proclaim{Proposition 3.2.2} ([\cite{10}, chapter I, section 1.1.4]) If
$X$ is a nonalgebraic surface and $\kappa(X)=0$. Then $\pi_1(X)$ lies in
an exact sequence. $$1\to {\Bbb Z}\times \pi_1(N)\to \pi_1(X)\to G\to 1$$
where $N$ is a $3$-dimensional nilmanifold and $G$ is a finite
group.\endproclaim

\subhead
3.3 The case $\kappa(X)=1$
\endsubhead

An {\it elliptic structure} or {\it elliptic fibration} ([\cite{10},
chapter I, section 1.14]) on a complex surface $X$ is a holomorphic map
$\pi:X\to C$, where $C$ is a curve such that for a general point $t\in C$,
$\pi^{-1}(t)$ is a curve of genus $1$. An {\it elliptic surface} is
a surface with an elliptic structure. Note that given an elliptic
fibration the base curve $C$ inherits an orbifold structure. 

\proclaim{Theorem 3.3.1} ([\cite{10}, chapter I, theorem 1.15]) If $X$ is
a complex surface with $\kappa(X)=1$ then $X$ is an elliptic
surface.\endproclaim

\subhead
3.4 The case $\kappa(X)=2$
\endsubhead

Surfaces with $\kappa(X)=2$ are called surface of {\it general type}.
There is no complete classification of these surfaces. As examples
of these surfaces we recall some of them from [\cite{2}, p. 189]. $(a)$
complete intersections of sufficiently high degree; $(b)$ products of
curves of genus $\geq 2$; $(c)$ Kodaira fibration; $(d)$ quotients of
symmetric domains; and $(e)$ practically all ramified double covering of
${\Bbb CP}^2$. In the last case the surface can be a singular algebraic
surface. 

\subhead
3.5 On elliptic surfaces
\endsubhead

We have already defined elliptic surfaces. In this subsection we
recall the structure of the fundamental group of an elliptic surface. Also
we recall the definition of $C^{\infty}$-elliptic surface and describe
their fundamental groups. We always follow the terminology of \cite{10}.

\proclaim{Proposition 3.5.1} ([\cite{10}, chapter II, theorem 2.3 and
proposition 2.1]) Let $\pi:X\to C$ be an elliptic fibration. If the Euler
number of $X$ is positive then $\pi$ induces an isomorphism $\pi_1(X)\to
\pi_1^{\text orb}(C)$.\endproclaim

\proclaim{Proposition 3.5.2} ([\cite{10}, chapter II, lemma 7.3 and
proposition 7.4]) Let $\pi:X\to C$ be an elliptic surface with Euler
number $0$. Then we have 

\roster
\item if $\pi_1^{\text orb}(C)$ is infinite then there is an exact
sequence $$0\to {\Bbb Z}\oplus {\Bbb Z}\to \pi_1(X)\to \pi_1^{\text
orb}(C)\to 1.$$
\item if $\pi_1^{\text orb}(C)$ is finite and the orbifold $C$ is good
then there is an exact sequence $${\Bbb Z}\to {\Bbb Z}\oplus {\Bbb Z}\to
\pi_1(X)\to \pi_1^{\text orb}(C)\to 1.$$
\item if $C$ is bad then $\pi_1(X)$ is abelian and is isomorphic to either
${\Bbb Z}\oplus {\Bbb Z}$ or ${\Bbb Z}\oplus {\Bbb Z}/{n{\Bbb Z}}$.
\endroster
\endproclaim

\proclaim{Definition 3.5.3} ([\cite{10}, chapter II, definition 1.1]) {\rm
A $C^{\infty}$-{\it elliptic 
surface} is a smooth map $\pi:X\to C$ from a closed smooth oriented
$4$-manifold $X$ to a smooth oriented real $2$-manifold $C$ such that for
each point $p\in C$ there is an open disk $p\in \Delta\subset C$ and a
complex elliptic surface $S\to \Delta$ and a smooth orientation preserving
diffeomorphism $\pi^{-1}(\Delta)\to S$ such that the following diagram
commutes.} 

$$\diagram
\pi^{-1}(\Delta)\rto^{\simeq} \dto & S \dto\\
\Delta \rto^{=} &\Delta\enddiagram$$
\endproclaim

It follows that for a $C^{\infty}$-elliptic surface $X$ the general fiber
of the map $\pi:X\to C$ is a $2$-torus. We say $X$ has no {\it singular
fiber} if all the fibers of $\pi:X\to C$ are smooth submanifold. Note that
$C$ inherits a orbifold structure and when the Euler number of $X$ is zero 
the monodromy action of $\pi_1^{\text orb}(C, p)$ on $H_1(f)$, where $f$
is the fiber of $\pi$ over a general base point $p$, factors through
$\pi_1(C)$ and its image is a finite cyclic subgroup of $SL(H_1(f))$. In
such a situation we say that the elliptic surface $X$ has {\it cyclic
monodromy} (see last paragraph in [\cite{10}, p. 195]). 

Now we are ready to state our next proposition whose proof is easily
deduced from the proof of Proposition 3.5.2.

\proclaim{Proposition 3.5.4} Let $X$ be a $C^{\infty}$-elliptic surface
without any singular fiber and has cyclic monodromy. Then the conclusions 
of Proposition 3.5.2 are true.\endproclaim

\head
4. Proofs
\endhead

\demo{Proof of Theorem 1.2} If $X$ is algebraic and $\kappa(X)=-\infty$
then by Proposition 3.1.1 either $\pi_1(X)$ is trivial or isomorphic to the
fundamental group of a curve. Hence the FIC is true in this case (Lemma
B). 

By Proposition 3.1.2 any Hopf surface has virtually infinite cyclic
fundamental group and hence the FIC is true for $\pi_1(X)$ by
[\cite{7}, lemma 2.7].

When $X$ is algebraic and Kodaira dimension $0$ then by Proposition 3.2.1
$\pi_1(X)$ is virtually abelian and hence the FIC is true
by [\cite{7}, lemma 2.7]. If it is nonalgebraic and $\kappa(X)=0$ then
by Proposition 3.2.2 $\pi_1(X)$ can be embedded as a cocompact discrete
subgroup of a virtually connected Lie group using [\cite{9}, algebraic
lemma]. In short we have $$\pi_1(X) < ({\Bbb Z}\times \pi_1(N))\wr
G < {\text Iso}(({\Bbb R}\times N)^G).$$

Hence the FIC is true for $\pi_1(X)$ by Theorem E. 

Inoue surfaces has poly-$\Bbb Z$ fundamental group by Corollary 3.1.4.
Hence by Theorem F the FIC is true for the fundamental group of an
Inoue surface.

Now we come to the case of elliptic surfaces. In the case when the Euler
number of the elliptic fibration is positive then by Proposition 3.5.1 
the fundamental group is isomorphic to the fundamental group of a
$2$-dimensional real orbifold and hence the FIC is true by Lemma
D. 

Next assume that the elliptic fibration has Euler number $0$. Then we had
three cases as in Proposition 3.5.2. By Lemma D we know that
$\pi_1^{\text orb}(C)$ satisfies the FIC. Assume $\pi_1^{\text orb}(C)$ is
infinite. Let $E$ be a virtually cyclic subgroup of $\pi_1^{\text
orb}(C)$. If $E$ is finite then $p^{-1}(E)$ is virtually abelian and hence
the FIC is true for $p^{-1}(E)$ ([\cite{7}, lemma 2.7]). Here $p$ is
the surjective homomorphism $p:\pi_1(X)\to \pi_1^{orb}(C)$. Now let
$Z$ be an infinite cyclic normal subgroup of $E$ of finite index. Then the
filtration $0<{\Bbb Z}<{\Bbb Z}\oplus {\Bbb
Z}<p^{-1}(Z)<p^{-1}(E)$ of $p^{-1}(E)$ gives a virtually poly-$\Bbb Z$  
group structure on $p^{-1}(E)$ and hence the FIC is true for $p^{-1}(E)$.
Thus
the FIC is true for $\pi_1(X)$. 

If $\pi_1^{\text orb}(C)$ is finite then from the second exact sequence in
Proposition 3.5.2 it follows that $\pi_1(X)$ is virtually abelian. Hence
the FIC is true for $\pi_1(X)$ by [\cite{7}, lemma 2.7]. In the last case
again we apply [\cite{7}, lemma 2.7]. 

Thus we have proved the Theorem 1.2.\qed\enddemo

\demo{Proof of Corollary 1.4} The proof is same as the proof of Theorem
1.2 in case when $X$ is an elliptic fibration.\qed\enddemo

\demo{Proof of Theorem 1.3} Recall that we have a fiber bundle projection 
from the $4$-dimensional real manifold $X$ over the circle and also $X$
has a complex structure. The following theorem shows that the fiber of
this fiber bundle projection is of a particular type. The theorem is a
consequence of [\cite{12}, theorem 7] and the discussion following it. 

\proclaim{Theorem 4.1} Let $X$ be a complex surface which is the
total space of a fiber bundle over ${\Bbb S}^1$. Then the fiber is
diffeomorphic to a Seifert fibered space.\endproclaim

Let $N$ be the fiber and $f$ be the monodromy diffeomorphism. 

It is well known (see [\cite{13}, theorem VI.17] and \cite{29}, \cite{30})
that apart from the following few cases ({\bf A} to {\bf E}) (see
[\cite{13}, section VI.16]) the diffeomorphism $f$ is isotopic to a fiber
preserving diffeomorphism, say $f$ again, that is $f$ sends the fiber
circle to fiber circle of $N$.

\noindent
{\bf A.} The Lens spaces including ${\Bbb S}^2\times {\Bbb S}^1$ and
${\Bbb S}^3$.

\noindent
{\bf B.} Seifert fibered spaces with base ${\Bbb S}^2$ and three
exceptional fiber of index corresponding to the triple $(2, 2, \alpha)$
where $\alpha>1$. These manifolds are called {\it prism-manifolds} and
has finite fundamental group.

\noindent
{\bf C.} A class of torus bundles over the circle. 

\noindent
{\bf D.} The solid torus. 

\noindent
{\bf E.} A twisted $I$-bundle over the Klein bottle.

At first let us assume that $f$ is fiber preserving. 
We have two cases.

{\bf $\pi_1(N)$ is infinite.} In this case 
as $N$ is a Seifert fibered space there is an exact sequence. $$1\to 
Z\to \pi_1(N)\to \pi_1^{orb}(B)\to 1$$ where $Z$ is infinite cyclic and is
generated by a regular fiber and $\pi_1^{orb}(B)$ is the orbifold
fundamental group of the base orbifold $B$ of $N$. 

As $f$ is fiber preserving and $Z$ is generated by a regular fiber the
induced action of $f$ on $\pi_1(N)$ leaves $Z$ invariant and hence we get
the following.
$$1\to Z \to \pi_1(N)\rtimes \langle t\rangle \to
\pi_1^{orb}(B)\rtimes \langle t\rangle \to 1.$$

Which reduces to: $$1\to Z \to \pi_1(X)\to 
\pi_1^{orb}(B)\rtimes \langle t\rangle \to 1.$$ Here, up to conjugation,  
the actions of $t$ on the various groups are induced by $f$. 

We will apply Lemma C to this exact sequence. If $\pi_1^{orb}(B)$ is
finite then $\pi_1(X)$ is virtually poly-$\Bbb Z$ and hence the FIC is
true for this group by Theorem F. So we assume $\pi_1^{orb}(B)$ is
infinite. In this case 
there is a finite index subgroup $H$ of $\pi_1^{orb}(B)$
which is the fundamental group of a closed surface. 

\noindent
{\bf Claim.} There is a characteristic closed surface subgroup of
$\pi_1^{orb}(B)$ of finite index.

\demo{Proof of claim} Let $H_1$ be the intersection of all conjugates of
$H$ in $\pi_1^{orb}(B)$. Then $H_1 < H$ is a finite index normal subgroup
of $\pi_1^{orb}(B)$. Also this implies that $H_1$ is again a closed
surface group. Now let $G=\pi_1^{orb}(B)/H_1$ be the finite quotient
group. Since $\pi_1^{orb}(B)$ is finitely presented there are only
finitely many homomorphism from $\pi_1^{orb}(B)$ to $G$, say $f_1,\cdots ,
f_n$. Consider $K=\cap_i\text{ker}f_i$. Then $K$ is a finite 
index characteristic subgroup of $\pi_1^{orb}(B)$. Also $K < H_1$ and
hence $K$ is also a closed surface group.\qed\enddemo  

\remark{Remark 4.2} {\rm We will use the procedure of finding a finite
index characteristic subgroup as in the proof of the above claim in the
remaining part of the paper without repeating the argument  
again.}\endremark

Since $K$ is a characteristic subgroup we have an exact sequence. $$1\to
K\to \pi_1^{orb}(B)\rtimes \langle t \rangle \to (\pi_1^{orb}(B)/K)\rtimes
\langle t \rangle \to 1.$$ 

Note that $(\pi_1^{orb}(B)/K)\rtimes \langle t \rangle $ is a
virtually cyclic group. Let $C$ be an infinite cyclic normal subgroup
of $(\pi_1^{orb}(B)/K)\rtimes \langle t \rangle $ of finite index. Then we
get the following. $$\pi_1^{orb}(B)\rtimes \langle t \rangle < (K\rtimes
C)\wr H.$$ Where $H\simeq ((\pi_1^{orb}(B)/K)\rtimes \langle t \rangle ) 
/C$ is a finite group. Since $K$ is the fundamental group of a closed
surface $K\rtimes C$ is the fundamental group of a closed $3$-manifold
which fibers over ${\Bbb S}^1$. If $K$ is the fundamental group of a
torus then $(K\rtimes C)\wr H$ is virtually poly-$\Bbb Z$ and hence the 
FIC is
true. Otherwise, by hypothesis there is a base diffeomorphism (the
monodromy diffeomorphism) which is special and therefore we can apply
Proposition 1.7 to deduce that the FIC is true for $(K\rtimes C)\wr H$.
Consequently, the same is true for $\pi_1^{orb}(B)\rtimes \langle
t\rangle$. To complete the proof we need to check that $p^{-1}(C)$
satisfies the FIC for any virtually cyclic subgroup $C$ of
$\pi_1^{orb}(B)\rtimes \langle t\rangle$. Here $p$ denotes the
homomorphism $\pi_1(X)\to \pi_1^{orb}(B)\rtimes \langle t\rangle$.
Obviously $p^{-1}(C)$ is virtually poly-$\Bbb Z$ and hence the FIC is true
for
$p^{-1}(C)$. This completes the proof in this case.

{\bf $\pi_1(N)$ is finite.} As $\pi_1(N)$ is finite $\pi_1(X)\simeq
\pi_1(N)\rtimes \langle t\rangle $ has a normal subgroup of finite index
isomorphic to
$\pi_1(N)\times {\Bbb Z}$. Hence $\pi_1(X)$ is virtually infinite cyclic.
The theorem now follows.

Now we come to the situation when $f$ is not fiber preserving. We have
mentioned before there are the following possibilities: {\bf A} to {\bf
E}. 

\noindent
{\bf Case A.} $\pi_1(X)$ is isomorphic to one of the following. 

\roster
\item fundamental group of a closed flat $2$-manifold.
\item infinite cyclic.
\item has a finite normal subgroup with infinite cyclic quotient. And
hence is virtually infinite cyclic.
\endroster 

We have already mentioned that the FIC is true in cases $(1)$ to $(3)$.

\noindent
{\bf Case B.} This case goes to $(3)$ in Case A.

\noindent
{\bf Case C.} In this case $\pi_1(X)$ is poly-$\Bbb Z$ and hence the FIC
is true by Theorem F.

As $X$ is closed Case D and Case E do not occur.

This completes the proof of the Theorem.\qed\enddemo

\remark{Remark 4.3} {\rm The crucial fact we used in the proof is that
the fiber $N$ is a Seifert fibered space. We would like to point out here
that except for only two examples all orientable Seifert fibered space
appear as fiber of fiber bundle projection $M^4\to {\Bbb S}^1$ where
$M^4$ is a complex surface (see [\cite{12}, theorem 7].)}\endremark

\demo{Proof of Corollary 1.8} If $G$ is torsion free and the FIC is true
for
$G$ then using Lemma C to the projection $G\times {\Bbb Z}^n\to G$ 
and noting that the FIC is true for free abelian groups we get that the 
FIC is
true for $G\times {\Bbb Z}^n$ also. Now we can apply [\cite{9}, theorem
D] to see that $Wh(G\times {\Bbb Z}^n)=0$ for all $n\geq 0$. It now
follows from Bass's contracted functor argument ([\cite{1}, \S 7, chapter
XII]) that $\tilde K_0({\Bbb Z}G)=K_{-i}({\Bbb Z}G)=0$ for all $i\geq
1$.\qed\enddemo 

In this paragraph we will say few words about the case of surfaces of
general type. Recall the examples $(a)$ to $(e)$ from 3.4 of Section 3.
The examples in $(a)$ are simply connected ([\cite{2}, chapter V,
proposition 2.1]). Examples $(b)$ and $(c)$ are already considered in
Theorem 1.5. The fundamental groups of compact quotients of symmetric
domains are discrete cocompact subgroups of the Lie groups $SU(2,1)$ or
$PSL(2, {\Bbb R})\times PSL(2, {\Bbb R})$ (see [\cite{2}, chapter V,
section 20]). In case $(e)$ we
have complex surfaces $X$ with a surjective morphism $p:X\to {\Bbb
{CP}}^2$ such that outside a curve $C$ in ${\Bbb {CP}}^2$ the morphism is
a $2$-sheeted covering; that is $p$ is a ramified covering with
ramification locus $C$. It is well known that for any curve $C'\subset   
X$ the inclusion map induces a surjective homomorphism $\pi_1(X-C')\to 
\pi_1(X)$. Now combining the two results from [\cite{21},
corollary(Zariski
Conjecture) and corollary 2.8]  we get that
if $C\subset {\Bbb {CP}}^2$ is either a smooth curve or a singular curve
with only nodal singularities then $\pi_1({\Bbb {CP}}^2-C)$ is abelian.
Hence $X-p^{-1}(C)$ is also abelian. On the other hand we have a
surjective homomorphism $\pi_1(X-p^{-1}(C))\to \pi_1(X)$. Thus we have
proved the following.

\proclaim{Proposition 4.4} Let $X$ be an algebraic surface (possibly
singular). Then $\pi_1(X)$ is abelian whenever $X$ is a $2$-sheeted
ramified covering of ${\Bbb {CP}}^2$ with ramification locus is either a
smooth or a nodal curve.\endproclaim

Summarizing the above discussions we note that for these surfaces of
general type the fundamental groups are equal to $1$, groups already
contained in Theorem 1.5, discrete cocompact subgroup of Lie group or
abelian. We have already shown that the FIC is true in these situations.

\remark{Remark 4.5} {\rm Recall from [\cite{2},chapter I, section 17]  
that the algebraic surface in the above proposition is singular whenever
the branch locus is singular.}\endremark

We give an example below to show that, given the method of proof we have
to prove the FIC for groups, among the examples of ramified coverings of
${\Bbb {CP}}^2$ the above class of examples are best possible for which
we could prove the FIC. 

\proclaim{Examples 4.6} ([\cite{21}, example 6.7]) {\rm Let $f$ and
$g$ be two homogeneous polynomials in three variables of degree $2$ and
$3$ respectively. Consider the curve $C$ defined by $f^3-g^2=0$ in ${\Bbb
{CP}}^2$. $C$ is smooth outside the points where $f=0$ and $g=0$. And the
singular points are ordinary cusp. Then $\pi_1({\Bbb {CP}}^2-C)\simeq
{\Bbb Z}_2 * {\Bbb Z}_3$. Now if $X$ is a complex surface which is a
$2$-sheeted ramified covering of ${\Bbb {CP}}^2$ with ramification locus
$C$ then $\pi_1(X)$ is the image of an index $2$ subgroup $H$ of ${\Bbb
Z}_2 * {\Bbb Z}_3$. On the other hand $H$ contains the commutator subgroup
of ${\Bbb Z}_2 * {\Bbb Z}_3$ which is a nonabelian free group. Hence
$\pi_1(X)$ can be a highly complicated group.}\endproclaim  

\head
5. The FIC for some real manifolds
\endhead

We need the following crucial Lemma to prove Theorem 1.5. Apart from
being the crucial ingredient throughout the paper the Main Lemma is of
independent interest. The proof of this Lemma is given in Section 6.
The definition of special diffeomorphism is given in the Introduction.

\proclaim{Main Lemma} Let $M^3$ be a closed $3$-dimensional 
manifold which is the total space of a fiber bundle projection $M^3\to
{\Bbb S}^1$ with orientable fiber. When the fiber is a surface of genus
$\geq 2$ assume that the fiber bundle has special monodromy
diffeomorphism. Then the FIC is true for $\pi_1(M)$.\endproclaim

\demo{Proof of Theorem 1.5} We prove a general version of Theorem 1.5. 
Let $F$ and $B$ be two closed orientable $2$-dimensional real manifold.
Let $G$ be a group which fits in an exact sequence. $$1\to \pi_1(F)\to
G\to \pi_1(B)\to 1.$$ Also by assumption the action of any element of
$\pi_1(B)$ on $\pi_1(F)$ is induced by a special diffeomorphism of $F$.

Then we prove that the FIC is true for $G$. Note that $\pi_1(X)$ also fits
in
such an exact sequence. Without loss of generality we can assume that
both the surfaces $F$ and $B$ are of genus $\geq 1$. Then as $B$ supports
a nonpositively curved Riemannian metric, $\pi_1(B)$ satisfies the FIC.
Let
$C$ be a virtually cyclic subgroup of $\pi_1(B)$. As $\pi_1(B)$ is torsion
free, $C$ is infinite cyclic. We have an exact sequence $1\to \pi_1(F)\to
p^{-1}(C)\to C\to 1$. We have $p^{-1}(C)\simeq \pi_1(F)\rtimes \langle
t\rangle $ where
the image of $t$ in $C$ generates $C$. It is well known that as $F$  is
closed, up to conjugation the action of $t$ on $\pi_1(F)$ is induced by a
diffeomorphism of $F$. Hence $p^{-1}(C)\simeq
\pi_1(N)$ where $N$ is a $3$-manifold which is a fiber bundle over the
circle with fiber diffeomorphic to $F$. By hypothesis this
diffeomorphism is special. By the Main Lemma the FIC is
true for $p^{-1}(C)$. We apply Lemma C to complete the proof.\qed\enddemo

\demo{Proof of Corollary 1.6} The proof is the same as the proof of
Theorem 1.5 when $M$ is compact. The only exception will be if the fiber
has nonempty boundary. If the fiber is simply connected then there is
nothing to prove. Otherwise by theorems 3.2 and 3.3 from \cite{18} the
interior of the mapping torus $M$ of a fiber supports a complete
nonpositively curved Riemannian metric so that near the boundary the
metric is a product metric. Hence the double of $M$ will support a
nonpositively curved metric. Also, the mapping torus $M$ will
have incompressible boundary and hence $\pi_1(M)<\pi_1(DM)$ where $DM$
denotes the double of $M$. By Lemma B we complete the proof in this case.

Next we assume $M$ is noncompact and of finite topological type. Consider
the surjective homomorphism $p:\pi_1(M)\to \pi_1(N)$ with kernel
$\pi_1(F)$. If $C$ is an infinite cyclic subgroup of $\pi_1(N)$ then
$p^{-1}(C)\simeq \pi_1(F)\rtimes {\Bbb Z}\simeq \pi_1(P)$ where $P$ is the
mapping torus of a (monodromy) diffeomorphism of $F$ . Note that
$\pi_1(F)$ is a finitely generated free group. Now we use Lemma 6.1 and
the previous argument in case of compact fiber with nonempty boundary to
complete the proof of the Corollary in this case. If the
fiber has infinitely generated fundamental group then
$\pi_1(P)$ is a subgroup of a $3$-manifold fibering over the
circle with special monodromy diffeomorphism and hence the FIC is true for
$\pi_1(P)$ by the Main Lemma and Lemma A.\qed\enddemo

\head
6. Proof of the Main Lemma
\endhead

\demo{Proof of the Main Lemma} We start the proof without assuming the
monodromy is special. The following exact sequence is obtained from the
long exact homotopy sequence of the fibration $M\to {\Bbb S}^1$. $$1\to
\pi_1(F)\to \pi_1(M)\to \pi_1({\Bbb S}^1)\to 1$$ where $F$ is the fiber of
the fiber bundle projection. Let $[A,A]$ denotes the commutator subgroup
of the group $A$. Then we have $$1\to [\pi_1(F), \pi_1(F)]\to \pi_1(F)\to
H_1(F, {\Bbb Z})\to 1.$$ Let $t$ be a generator of $\pi_1({\Bbb S}^1)$.
Since $[\pi_1(F), \pi_1(F)]$ is a characteristic subgroup of $\pi_1(F)$
the action (induced by the monodromy) of $t$ on $\pi_1(F)$ leaves
$[\pi_1(F), \pi_1(F)]$ invariant. Thus we have another exact sequence 
$$1\to [\pi_1(F), \pi_1(F)]\to \pi_1(F)\rtimes
\langle t\rangle \to H_1(F, {\Bbb Z})\rtimes \langle t\rangle \to 1.$$
Which reduces to the sequence $$1\to [\pi_1(F), \pi_1(F)]\to \pi_1(M)\to
H_1(F, {\Bbb Z})\rtimes \langle t\rangle \to 1.$$

We would like to apply Lemma C to this exact sequence. 

If the fiber is ${\Bbb S}^2$ or ${\Bbb T}^2$ then $\pi_1(M)$ is
poly-$\Bbb Z$ and hence the FIC is true for $\pi_1(M)$. So assume that
the fiber has genus $\geq 2$.

Note that the group $H_1(F, {\Bbb Z})\rtimes \langle t\rangle $ is
poly-$\Bbb Z$. Hence by Theorem F the FIC is true for $H_1(F,
{\Bbb Z})\rtimes \langle t\rangle$. Let $C$ be a virtually cyclic subgroup
of $H_1(F, {\Bbb Z})\rtimes \langle t\rangle$. Let $p$ denotes
the surjective homomorphism $\pi_1(M)\to H_1(F, {\Bbb Z})\rtimes \langle
t\rangle$. We will show that the FIC is true for $p^{-1}(C)$. Note that
$C$ is either trivial or infinite cyclic. 

{\bf Case $C=1$}. In this case we have that $p^{-1}(C)$ is 
a nonabelian free group and hence is the fundamental group of a surface.
We need the following lemma to complete the proof of this case.

\proclaim{Lemma 6.1} Let $\Gamma$ be the fundamental group of a surface
then the FIC is true for $\Gamma\wr G$ for any finite group
$G$.\endproclaim 

\demo{Proof} If $\Gamma$ is finitely generated then $\Gamma$ is the
fundamental group of a compact surface and hence
the lemma follows from Lemma B. In the infinitely generated
case $\Gamma\wr G\simeq \lim_{i\to \infty} (\Gamma_i\wr G)$ where each
$\Gamma_i$ is a finitely generated free group. Now recall that any
finitely generated nonabelian free group is the fundamental group of a
compact surface with nonempty boundary. By Lemma B, and [\cite{7},
theorem 7.1] the proof is complete.\qed\enddemo

{\bf Case $C\neq 1$}. We have $p^{-1}(C)=[\pi_1(F),
\pi_1(F)]\rtimes \langle s\rangle $ where $s$ is a generator of $C$. Let
$\tilde F$ be the covering space of $F$ corresponding to the commutator
subgroup $[\pi_1(F), \pi_1(F)]$. 

Note that the monodromy diffeomorphism of $F$ lifts to a   
diffeomorphism of $\tilde F$ which in turn, up to conjugation, induces the
action of $t$ on $[\pi_1(F), \pi_1(F)]$. Also, the
induced action of $t$ on $H_1(F, {\Bbb Z})$ is given by $t(v)=\tilde
f\circ v\circ {\tilde f}^{-1}$, where $\tilde f:\tilde F\to \tilde F$ is
a lift of the monodromy diffeomorphism and $v\in H_1(F, {\Bbb Z})$. Here
$H_1(F, {\Bbb Z})$ is identified with the group of covering transformation
of the covering $\tilde F\to F$. From this observation it
follows that, up to conjugation, the action of $s$ on $[\pi_1(F),
\pi_1(F)]$ is induced by a diffeomorphism (say $f$) of $\tilde F$.
Indeed, if $s=(s_1, t^k)\in H_1(F, {\Bbb Z})\rtimes \langle t\rangle$
then $f=s_1\circ {\tilde f}^k:\tilde F\to \tilde F$.

Let $g:F\to F$ be a diffeomorphism of $F$ so that $f$ and a lift of
$g$ induce the same outer automorphism of $\pi_1(\tilde F)$.
Such a diffeomorphism exists. Namely, let $u$ be an element of
$\pi_1(F)\rtimes \langle t \rangle $ which goes to $s$. Then the
conjugation action by $u$ on $\pi_1(F)\rtimes \langle t \rangle $ leaves
$\pi_1(F)$ invariant. Since $F$ is closed there is a diffeomorphism $g$ of
$F$ which, up to conjugation, induces this action of $u$ on $\pi_1(F)$.
Clearly $g$ has the required property. Let $M_g$ and $M_f$ be the mapping
tori of $g$ and $f$ respectively. Then $\pi_1(M_f)$ is a subgroup of
$\pi_1(M_g)$ and $\pi_1(M_f)\simeq [\pi_1(F), \pi_1(F)]\rtimes \langle
s\rangle $.

Note that the topological type of the mapping torus $M_g$ is an invariant
of the isotopy class of $g$. Below, whenever we say ``$g$ is'' we mean
``an isotopy of $g$ is''. Also assume that $g$ is an orientation
preserving diffeomorphism. The orientation reversing case can easily be
tackled from the orientation preserving case and is left to the reader. 

By Nielsen-Thurston classification of surface diffeomorphisms
(see \cite{3} or [\cite{24}, p. 175]) there are now three cases. 

\noindent
{\bf Case 1.} $g$ is pseudo-Anosov. In this case $M_g$ supports a
hyperbolic structure (\cite{22}, \cite{28}) and hence the FIC is true for
$\pi_1(M_g)$ by Lemma B. Thus the FIC is true for $\pi_1(M_f)$ also. 

\noindent
{\bf Case 2.} $g$ is of finite order. In this case $M_g$ has a
regular finite sheeted covering diffeomorphic to $F\times {\Bbb S}^1$.
Hence $\pi_1(M_g)$ is a subgroup of $\pi_1(F\times {\Bbb S}^1)\wr G_1$
where $G_1$ is a finite group. Since $F\times {\Bbb S}^1$ supports a
nonpositively curved Riemannian metric Lemma B applies again.

\noindent
{\bf Case 3.} $g$ is not pseudo-Anosov and of infinite order. By the above
mentioned classification of surface diffeomorphism, in this case $g$ is
reducible (see [\cite{3}, p. 75] for definition) and hence there are
finitely many nontrivial elements $h_1, h_2,\cdots , h_n$ in $\pi_1(F)$
represented by pairwise disjoint mutually nonparallel simple closed curves
(say, $C_1, C_2,\ldots , C_n$ respectively) on $F$ so that $g(\bigcup_i 
C_i)=\bigcup_i C_i$. Also there are pairwise disjoint tubular
neighborhoods $N(C_i)$ of $C_i$ and submanifolds $F_j$ for $j=1,2,\cdots,
l$ of $F$ satisfying the followings (see [\cite{8}, p. 219,
Th\'{e}or\`{e}me 4.2]). 

\roster
\item $F-\bigcup_{i=1}^nN(C_i)=\bigcup_{j=1}^lF_j$;
\item $g(F_j)=F_j$ for each $j=1,2,\cdots ,l$;
\item $g|_{F_j}$ is isotopic either to a pseudo-Anosov diffeomorphism or
to a finite order diffeomorphism.
\endroster

Here some $F_j$ may have more than one connected component. Choose a
positive integer $N$ so that $g^N(L)=L$ and $g^N|_L$ is isotopic
either to a pseudo-Anosov diffeomorphism or to a finite order 
diffeomorphism, for each connected component $L$ of
$F-\bigcup_{i=1}^nN(C_i)$. Hence the mapping torus of $g^N|_L$ is either a
Seifert fibered space ( in the case when $g^N|_L$ is isotopic to a finite
order diffeomorphism) or supports a hyperbolic metric in the interior
(pseudo-Anosov case) [see \cite{28} or \cite{22}]. In fact the mapping
tori
of $g^N|_L$ where $L$ varies over the connected components of
$F-\bigcup_{i=1}^nN(C_i)$ give the JSJT decomposition of the mapping torus
$M_{g^N}$ of $g^N$. If for
some $L$, the mapping torus of $g^N|_L$ is hyperbolic then by [\cite{18},
theorem 3.2 and 3.3] $M_{g^N}$ supports a nonpositively curved Riemannian
metric and since $M_{g^N}$ is a finite sheeted cover of $M_g$ we can apply
Lemma B to conclude that the FIC is true for $\pi_1(M_g)$ and hence for
$\pi_1(M_f)$ (Lemma A) also. Therefore for the rest of the proof we can
assume that for each component $L$ of $F-\bigcup_{i=1}^nN(C_i)$ $g^N|_L$
is isotopic to a finite order diffeomorphism. From the above discussion it
follows that the FIC is true for $\pi_1(M_f)$ when the monodromy is
special and satisfies condition $(1)$ or $(2)$ in the definition. 

From now onwards we assume that the monodromy diffeomorphism is
special and satisfies condition $(3)$ in the definition. 

We will write $M_{f^{N}}$ as an increasing union of connected compact
$3$-manifolds with incompressible tori boundary components.

We need the following claim.

\noindent
{\bf Claim.} $M_f$ has one topological end.

\demo{Proof of claim} Recall that ([\cite{5}, p. 115]) by definition a
group $G$ has
$e(G)$ number of ends if there is a regular covering projection $\tilde
X\to X$ where $X$ is a finite complex, $G$ is isomorphic to the group of
covering transformation of $\tilde X\to X$ and $\tilde X$ has $e(G)$
number of topological ends. Also this definition does not depend on the
covering projection ([\cite{5}, theorem 3]).

As $F$ has first Betti number $\geq 2$ the group $H_1(F, {\Bbb Z})$ is
free abelian of rank greater than $1$ and also $H_1(F, {\Bbb Z})$ is
the group of covering transformation of ${\Bbb R}^n\to ({\Bbb S}^1)^n$
where $n$ is the rank of $H_1(F, {\Bbb Z})$. Hence $H_1(F, {\Bbb Z})$ has
one end. Also $H_1(F, {\Bbb Z})$ is the group of covering transformations
of the regular covering $\tilde F\to F$. Since $F$ is compact, the
manifold $\tilde F$ has one topological end (see \cite{5}). Figure 1
describes $\tilde F$.

\medskip

\centerline{\psfig{figure=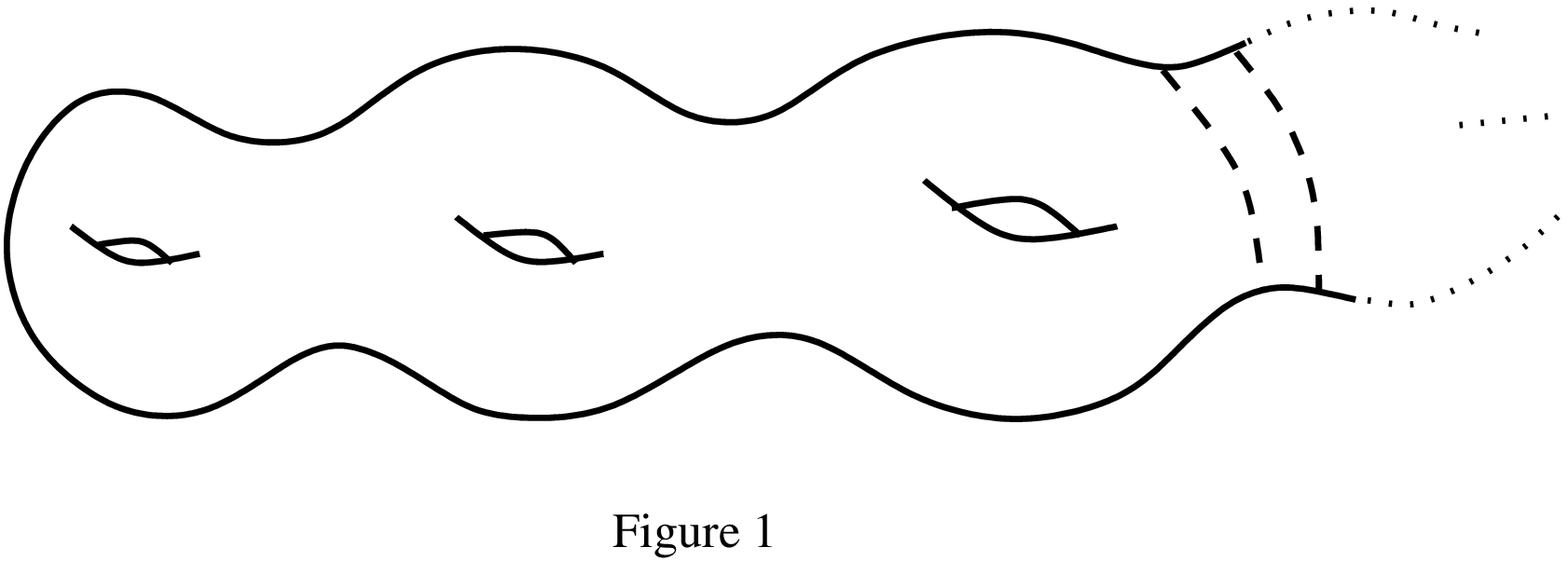,height=3cm,width=7cm}}

\vskip 0.7cm

It now follows that $M_f$ has one topological end.\qed\enddemo

The same proof shows that $M_{f^{N}}$ also has one topological end. This
ensures that there exists a surjective proper smooth function 
$\delta:M_{f^{N}}\to {\Bbb R}^{\geq 0}$. Also since $M_{f^{N}}$ is
connected and has one end we can always choose an $r\in {\Bbb R}^{\geq
0}$ bigger than any given positive real number such that 
$\delta^{-1}([0,r])$ is connected.

We would like to construct submanifolds $\tilde N_i^s$ of $M_{f^{N}}$ 
with the following properties. 

\roster
\item each $\tilde N_i^s$ is a compact connected submanifold of
$M_{f^{N}}$ and has incompressible tori boundary.
\item $\tilde N_i^s\subset \tilde N_{i+1}^s$ for $i=1,2, \cdots$.
\item $M_{f^{N}}=\cup_{i=0}^{i=\infty}\tilde N_i^s$.
\endroster

Consider the covering projection $p:M_{f^{N}}\to M_{g^{N}}$. Recall that
all the pieces in the JSJT decomposition of $M_{g^{N}}$ are Seifert
fibered and are obtained by taking the mapping torus of $g^{N}|_L$ for
some component $L$ of $F-\bigcup_{i=1}^nN(C_i)$. Let $P$ be such a Seifert
fibered piece. We claim the following. 

\noindent
{\bf Claim.} each component of $p^{-1}(P)$ is a Seifert fibered space. 

\demo{Proof of claim} Recall that $u\in \pi_1(F)\rtimes \langle
t\rangle$ was a lift of $s=(s_1, t^k)\in H_1(F, {\Bbb Z})\rtimes \langle 
t\rangle$ and $g$ induced the restriction to $\pi_1(F)$ of the conjugation
action by $u$ on $\pi_1(F)\rtimes \langle t\rangle$. Let $u=(u_1, t^k)$
and let $(x, 1)\in \pi_1(F)\rtimes \langle t\rangle$. Then
$u(x,1)u^{-1}=(u_1t^k(x)u_1^{-1}, 1)$. That is, the action of $u$ is
composition of the action of $t^k$ and a conjugation by an element of
$\pi_1(F)$. Since conjugation action by an element of $\pi_1(F)$ is
induced by a diffeomorphism of $F$ isotopic to the identity  
(\cite{20}), it follows that $g$ is isotopic to the $k$-th power of the
monodromy diffeomorphism (say $h$). Hence $g^N$ is isotopic to
$h^{kN}$. Consequently $M_{h^{kN}}$ is diffeomorphic to $M_{g^N}$.

Consider the following commutative diagram. 

$$\diagram
M_{f^N} \rto^{f_1} \dto^p & M_{\tilde h} \dto^q\\
M_{g^N} \rto^{f_2} &M_h\enddiagram$$

As $M_{g^N}\to M_h$ is a finite sheeted cover the JSJ decomposition of
$M_{g^N}$ can be obtained by taking the inverse image of the pieces of the
JSJT decomposition of $M_h$. So let $P$ be such a piece in $M_{g^N}$
and let $S_P$ be a component of $p^{-1}(P)$ (as above) which goes to the
component $S_{P'}$ of $q^{-1}(P')$ where $P'=f_2(P)$. Then we have a
finite sheeted covering projection $S_P\to S_{P'}$ and by hypothesis
$\pi_1(S_{P'})$ is not free. We check that $S_{P'}$ is a Seifert fibered
space. Since $S_P\to S_{P'}$ is finite sheeted it will follow that $S_P$
is also Seifert fibered.

The claim now follows from the following Lemma.

\proclaim{Lemma 6.2} In the above notation $S_{P'}$ is a Seifert
fibered space.\endproclaim

\demo{Proof} Recall the well-known theorem that if the
fundamental group of a compact orientable irreducible $3$-manifold
contains an infinite cyclic central subgroup then the manifold admits a
Seifert fibered structure. There is a generalization of this theorem
for noncompact $3$-manifolds (see [\cite{19}, theorem 1.1]) under certain
conditions ([\cite{19}, definitions in the introduction]). Since $P'$
is a compact quotient of $S_{P'}$ these conditions are easily satisfied by
$S_{P'}$. Hence we only have to show that there is an infinite cyclic
central subgroup of $\pi_1(S_{P'})$. Since $P'$ is Seifert fibered we
have the following short exact sequence. $$1\to C\to \pi_1(P')\to
\pi_1^{orb}(B)\to 1$$ where $C$ is infinite cyclic central (generated by
a regular fiber) and $B$ is the base surface of the Seifert fibered
space $P'$. If $C\cap p_*(\pi_1(S_{P'}))=(1)$ then since $B$ has
nonempty boundary and $\pi_1(S_{P'})$ is torsion free it would follow that
$\pi_1(S_{P'})$ is free. Which is a contradiction since $h$ is special
and satisfies condition $(3)$.\qed\enddemo 

This proves the claim.\qed\enddemo

Thus we have found a JSJT type decomposition of $M_{f^{N}}$ consisting of
Seifert fibered spaces which are components of $p^{-1}(P)$ where
$P$ varies over all Seifert fibered pieces of $M_{g^{N}}$. Each
component $S_P$ of $p^{-1}(P)$ has base surface (say $B_S$) with a
discrete set (say, $D$ ) of orbifold points on $B_S$. Let $q:S_P\to
B_S$ be the quotient map. Note that $B_S$ has a
filtration by increasing sequence (under inclusion) of compact subsurfaces
with incompressible circle boundary components so that each such circle
boundary avoids $D$. By taking the inverse images of these subsurfaces
under $q$ we can write $S_P$ as an increasing union of 
compact submanifolds $N_i^{S_P}$ with incompressible tori boundary
components. 

Now choose $r_1\in {\Bbb R}^{\geq 0}$ and $i^1_1,i^1_2,\cdots , i^1_{k_1}$
so that $$\delta^{-1}([0, r_1])\subset \cup_{j=1}^{j=k_1} N_{i^1_j}^{S_
{P_j}}$$ where $P_j$ is some Seifert fibered piece of $M_{g^{N}}$ and
$S_{P_j}$ denotes some component of $p^{-1}(P_j)$.
Since $M_{f^{N}}$ is connected and the family $\{N_i^{S_P}\}$ is a
covering of $M_{f^{N}}$ we can choose $i^1_1,i^1_2,\cdots ,
i^1_{k_1}$ so that $\cup_{j=1}^{j=k_1} N_{i^1_j}^{S_{P_j}}$ is connected.
Let $N_1^s=\cup_{j=1}^{j=k_1} N_{i^1_j}^{S_{P_j}}$. Similarly we can
choose $r_1<r_2\cdots <r_m<\cdots$ and $i^l_1,i^l_2,\cdots , i^l_{k_l}$,
$l=1,2,\cdots$, so that $r_m\to\infty$ and $$\cup_{j=1}^{j=k_{l-1}}
N_{i^{l-1}_j}^{S_{P_j}}\subset \delta^{-1}([0, r_l])\subset
\cup_{j=1}^{j=k_l} N_{i^l_j}^{S_{P_j}}$$ and $\cup_{j=1}^{j=k_l}
N_{i^l_j}^{S_{P_j}}$ is connected. Write $N_l^s=\cup_{j=1}^{j=k_l}
N_{i^l_j}^{S_{P_j}}$. It follows that $N_l^s$ satisfies the properties
$(2)$ and $(3)$. To check $(1)$ note that $N_l^s$ is a connected and
compact $3$-manifold. Also any two $N_{i^l_j}^{S_{P_j}}$ either do not
intersect or intersect along incompressible tori boundary components or
along some incompressible annuli on some incompressible tori boundary
components. In any case it follows that $N_l^s$ has incompressible tori
boundary. This proves $(1)$.

Now we check that the FIC is true for $\pi_1(M_f)$.
By [\cite{18}, theorems 3.2 and 3.3] the interior of $N_i^s$ supports
a complete nonpositively curved Riemannian metric so that near the
boundary (tori) the metric is a product. Since the inclusion $N_i^s
\subset N_i^s\cup_{\partial} N_i^s$ have the obvious retractions we get
the following inclusion on fundamental groups: $\pi_1(N_i^s) < 
\pi_1(N_i^s\cup_{\partial} N_i^s)$ where $N_i^s\cup_{\partial} N_i^s$ is
the double of $N_i^s$. Also $N_i^s\cup_{\partial} N_i^s$ is a closed
nonpositively curved manifold (by \cite{18}). Hence the FIC is true for
$\pi_1(N_i^s)\wr G$ by Lemma A and B, where $G$ is a finite group. Also we
have $\pi_1(M_{f^{N}}) \simeq \lim_{i\to \infty} \pi_1(N_i^s)$ and hence
we get that $\pi_1(M_f)<\pi_1(M_{f^{N}})\wr G  <\lim_{i\to \infty}
(\pi_1(N_i^s)\wr G)$, where $G$ is the group of covering transformation of
$M_{f^{N}}\to M_f$.  And hence by [\cite{7}, theorem 7.1] the FIC is true
for $\pi_1(M_f)\simeq [\pi_1(F), \pi_1(F)]\rtimes \langle u\rangle $.

This completes the proof of the Main Lemma.\qed\enddemo

Lastly we record, for later application, the following corollaries which
are consequences of the proof of the Main Lemma. 

\proclaim{Corollary 6.3} Let $F$ be a closed orientable surface of
genus $ > 1$ and $g:F\to F$ be an orientation preserving reducible and
infinite order diffeomorphism. Also assume that all the pieces in the
JSJT decomposition of $M_{g^N}$ (induced by $g$) for some large $N$ are
Seifert fibered. Assume that $g$ is special. Let $f:\tilde F\to
\tilde F$ be a lift of the diffeomorphism $g$ under the covering $\tilde
F\to F$ corresponding to the commutator subgroup of $\pi_1(F)$. Then
$\pi_1(M_f)\simeq \pi_1(\tilde F)\rtimes \langle s \rangle$ is a subgroup
of $\lim_{i\to \infty}(\pi_1(N_i^s)\wr G)$. Where $M_f$ is the mapping
torus of $f$, the action of $s$ on $\pi_1(\tilde F)$ is induced by $f$,
$N_i^s$ are increasing sequence (under inclusion) of compact connected
irreducible $3$-manifold with incompressible tori boundary components and
$G$ is a finite group. 
\endproclaim 

\proclaim{Corollary 6.4} In the hypothesis of the above corollary if $g$
is either pseudo-Anosov or is a finite order diffeomorphism or there
is a hyperbolic piece in the JSJT decomposition of $M_{g^N}$ for some
large $N$ then $\pi_1(M_f)$ is a subgroup of $\pi_1(N)\wr G$ where $N$ is
a closed nonpositively curved Riemannian $3$-manifold and $G$ is a finite
group.\endproclaim

\remark{Remark 6.5} We will use the Main Lemma and the method of its'
proof to deduce the FIC for a large class of $3$-manifold groups in
\cite{27}.\endremark

\remark{Remark 6.6} Here we remark that when the diffeomorphism $g$
belongs to {\bf Case 3} then $g$ is isotopic to the $k$-th power of
the monodromy diffeomorphism of the fiber bundle projection $M\to {\Bbb 
S}^1$.\endremark

\head
7. Virtually fibered $3$-manifold and the FIC
\endhead

In this section we prove Proposition 1.7.

\demo{Proof of Proposition 1.7} Note that it is enough to prove that the
FIC is true for $\pi_1(M)\wr G$ for any finite group $G$.

Let $F$ be the fiber of the fiber bundle projection $M\to {\Bbb S}^1$.
Passing to a finite cover we can make sure $F$ is orientable and
the monodromy diffeomorphism of the fiber bundle projection $M\to {\Bbb
S}^1$ is orientation preserving. There are two cases to consider now.

\noindent
{\bf Case A.} Assume that the monodromy diffeomorphism belongs to the
first two classes in the definition of special diffeomorphism. In this
case by Corollary 6.4 and [\cite{26}, lemma B] the proposition
follows.  

\noindent
{\bf Case B.} Assume that the monodromy diffeomorphism is special and
condition $(3)$ is satisfied. 

We have the following exact sequence. $$1\to [\pi_1(F), \pi_1(F)]\to
\pi_1(M)\to H_1(F, {\Bbb Z})\rtimes \langle t\rangle \to 1.$$ 

If $F$ is the $2$-sphere or the torus then $\pi_1(M)\wr G$ is
virtually poly-$\Bbb Z$ and hence the FIC is true by Theorem F. So
assume $F$ is not the $2$-sphere or the torus. That is the genus of $F$ is
$\geq 2$.

Taking wreath product with $G$ the above exact sequence gives the
following. $$1\to ([\pi_1(F), \pi_1(F)])^G\to \pi_1(M)\wr G\to (H_1(F,
{\Bbb Z})\rtimes \langle t\rangle)\wr G\to 1.$$ 

Note that $(H_1(F, {\Bbb Z})\rtimes \langle t\rangle)\wr G$ is
virtually
poly-$\Bbb Z$ and hence the FIC is true for $(H_1(F, {\Bbb Z})\rtimes
\langle
t\rangle)\wr G$. 

By the following Lemma, [\cite{7}, theorem 7.1] and Lemma 6.1 it follows
that the FIC is true for $([\pi_1(F), \pi_1(F)])^G$.

\proclaim{Lemma 7.2} Let $G_1$ and $G_2$ be two groups and assume
the FIC is true for both $G_1$ and $G_2$ then the FIC is true for the
product
$G_1\times G_2$.\endproclaim

\demo{Proof}  Consider the projection $p_1:G_1\times G_2\to
G_1$. By [\cite{26}, lemma C] we need to check that the FIC is true for
$p_1^{-1}(C)$ for any virtually cyclic subgroup $C$ of $G_1$. Note that
$p_1^{-1}(C)=C\times G_2$. Now consider the projection $p_2: C\times
G_2\to G_2$. Again we apply [\cite{26}, lemma C]. That is we need to show
that the FIC is true for $p_2^{-1}(C')$ for any virtually cyclic subgroup
$C'$
of $G_2$. But $p_2^{-1}(C')=C\times C'$ which is virtually poly-$\Bbb 
Z$ and hence the FIC is true for $p_2^{-1}(C')$. This completes the proof
of
the lemma.\qed\enddemo

Let $Z$ be a virtually cyclic subgroup of $(H_1(F, {\Bbb Z})\rtimes
\langle t\rangle)\wr G$. If $Z$ is finite then $$p^{-1}(Z) < 
([\pi_1(F), \pi_1(F)])^G\wr Z < ([\pi_1(F),
\pi_1(F)])\wr (G\times Z).$$

Here $p$ is the surjective homomorphism $\pi_1(M)\wr G\to (H_1(F, {\Bbb
Z})\rtimes \langle t\rangle)\wr G$. Now Lemma 6.1 applies on
the right hand side group to show that the FIC is true for $p^{-1}(Z)$.

If $Z$ is infinite then let $Z_1$ be the intersection of
$Z$ with the torsion free part $(H_1(F, {\Bbb
Z})\rtimes \langle t\rangle)^G$. Hence $Z_1 \ \simeq \ \langle u\rangle$
is an infinite cyclic normal subgroup of $Z$ of finite index. We get

$$p^{-1}(Z)<(p^{-1}(Z_1))\wr
{Z/Z_1} \simeq (([\pi_1(F),
\pi_1(F)])^G\rtimes \langle u\rangle )\wr Z/Z_1 \simeq$$$$ (([\pi_1(F),
\pi_1(F)]\times [\pi_1(F), \pi_1(F)]\cdots\times [\pi_1(F), 
\pi_1(F)])\rtimes \langle u\rangle )\wr
Z/Z_1 = H\text{(say)}.$$ 

In the above display there are $|G|$ number of factors of $[\pi_1(F),
\pi_1(F)]$. Note that the action of $u$ on $([\pi_1(F), \pi_1(F)])^G$ 
is factorwise, that is, the $i$-th coordinate of $u$ acts on the $i$-th
factor of $([\pi_1(F), \pi_1(F)])^G$. Let $u=(u_1,\cdots ,u_{|G|})$.
Recall that if $u_j\neq 1$ then the action of $u_j$ on $[\pi_1(F),
\pi_1(F)]$ is induced by
the lift of a diffeomorphism (say $g_j$) of $F$. Without loss of
generality we can assume that $g_1,\cdots, g_k$ are diffeomorphism of type
as in Corollary 6.3  and $g_{k+1},\cdots, g_l$ are of type as in
Corollary 6.4 and the remaining $u_j$ are trivial element. Note
that here it follows from Remark 6.6 that for $j=1,2,\cdots, k$ each
$g_j$ has the same property as $g$. Applying
Corollaries 6.3 and 6.4  we get $$H < ((\lim_{i \to \infty
}((\pi_1(N^{u_1}_i)\wr  
G_1)\times\cdots\times(\pi_1(N^{u_k}_i)\wr
G_k))\times(\pi_1(N_{k+1})\wr G_{k+1})\times\cdots\times$$$$ 
(\pi_1(N_l)\wr
G_l)\times ([\pi_1(F), \pi_1(F)])^{|G|-l})\wr Z/Z_1 < (\lim_{i \to
\infty }((\pi_1(N^{u_1}_i)\wr
G_1)\times\cdots\times$$$$(\pi_1(N^{u_k}_i)\wr
G_k))\times (\pi_1(N_{k+1})\wr G_{k+1})\times\cdots\times
(\pi_1(N_l)\wr
G_l)\times ([\pi_1(F), \pi_1(F)])^{|G|-l})\wr Z/Z_1 $$$$< (\lim_{i \to
\infty
}((\pi_1(N^{u_1}_i)^{G_1}\times\cdots\times\pi_1(N^{u_k}_i)^{G_k}\times
\pi_1(N_{k+1})^{G_{k+1}}\times\cdots\times\pi_1(N_l)^{G_l})\wr
(G_1\times\cdots$$$$\cdots\times G_l))\times ([\pi_1(F),
\pi_1(F)])^{|G|-l})\wr
Z/Z_1$$ where $N^{u_j}_i$ are
irreducible $3$-manifolds with incompressible boundary as appeared in
Corollary 6.3, $N_j$ are closed nonpositively curved Riemannian
manifolds and $G_1,\cdots, G_l$ are finite groups. Let us denote the
group inside the limit of the last expression by $K_i$. Then the last line
becomes $$\lim_{i \to \infty} (K_i\wr (G_1\times\cdots\times
G_l)\times ([\pi_1(F), \pi_1(F)])^{|G|-l})\wr
Z/Z_1$$$$  < \lim_{i \to  \infty} (K_i^{G_1\times\cdots\times
G_l\times Z/Z_1}\wr
((G_1\times\cdots\times G_l)\wr Z/Z_1))\times ([\pi_1(F),
\pi_1(F)])^{|G|-l}\wr Z/Z_1.$$ Let $M_i^{u_j}$ be the
double of $N_i^{u_j}$. Then by \cite{Le} $M_i^{u_j}$ is a closed
nonpositively curved Riemannian manifold. Using the obvious retractions
we find that $\pi_1(N_i^{u_j})$ is a subgroup of $\pi_1(M_i^{u_j})$.    
Hence $K_i$ is a subgroup of the fundamental group of a
closed nonpositively curved Riemannian manifold, namely
$(M^{u_1}_i)^{G_1}\times\cdots\times (M^{u_k}_i)^{G_k}\times
N_{k+1}^{G_{k+1}}\times\cdots\times N_l^{G_l}$. And the last
inclusion follows from the following easy to verify Lemma. 

\proclaim{Lemma 7.3} Let $A$ and $B$ be two finite groups and $G$ is any
group, then $(G\wr A)\wr B$ is a subgroup of $G^{A\times B}\wr (A\wr
B)$\endproclaim

Now using [\cite{7}, theorem 7.1], Lemma 6.1 , Lemma 7.2 and the
following Lemma we complete the proof of the Proposition.\qed\enddemo

\proclaim{Lemma 7.4} Let $A$ and $B$ be two groups and $G$ be a finite
group then $(A\times B)\wr G$ is a subgroup of $(A\wr G)\times(B\wr
G)$.\endproclaim

\demo{Proof} The proof is easy and left to the reader.\qed\enddemo

\head
8. Examples
\endhead
This section is devoted in giving examples of special diffeomorphisms. 

\proclaim{Lemma 8.1} Let $f:F\to F$ be a special diffeomorphism. Then for
any positive integer $n$, $f^n$ is also a special 
diffeomorphism. Conversely, if $f^n$ is special for some $n$
then so is $f$.\endproclaim

\demo{Proof} If $f$ satisfies condition $(2)$ in the definition of 
special diffeomorphism then obviously $f^n$ is special and satisfies
condition $(2)$. Converse is also trivial. If $f$ satisfies condition
$(1)$ then since $M_{f^n}$ is a finite sheeted cover of $M_f$ it follows
$M_{f^n}$ also satisfies condition $(1)$. In fact the pull back metric
does the job. The converse direction follows from [\cite{17}, corollary
2.5]. 

So assume that $f$ satisfies condition $(3)$.

Consider the following commutative diagram.

$$\diagram
M_{\tilde {f^n}} \rto^{f_1} \dto^p & M_{\tilde f} \dto^q\\
M_{f^n} \rto^{f_2} &M_f\enddiagram$$

Note that Seifert fibered pieces of $M_{f^n}$ are inverse images of
Seifert fibered pieces of $M_f$ under the map $f_2$. So let $P$ be a
Seifert fibered piece of $M_{f^n}$ and $f_2(P)=Q$. Let $\tilde P$ be a
component of $p^{-1}(P)$ which goes to the component $\tilde Q$ of
$q^{-1}(Q)$. Note that $f_1|_{\tilde P}:\tilde P\to \tilde Q$ is a finite
sheeted covering. By hypothesis $\pi_1(\tilde Q)$ is not free. By Lemma
6.2 $\tilde Q$ is Seifert fibered. Since $\tilde P$ is a finite sheeted
cover of $\tilde Q$, $\tilde P$ is also Seifert fibered and hence has a
infinite cyclic normal subgroup generated by a regular fiber. Hence
$\pi_1(\tilde P)$ is not free. Hence $f^n$ is also special. 

Conversely if $f^n$ is special for some $n$ then in the above
notation $\pi_1(\tilde P)$ is not free, consequently $\pi_1(\tilde Q)$ is
also not free as $\pi_1(\tilde P)$ is a subgroup of $\pi_1(\tilde Q)$. 
\qed\enddemo

\remark{Example 8.2} Now we give examples of some special diffeomorphisms.
We have already seen examples of special diffeomorphism which satisfies
condition $(1)$ or $(2)$. Let $f:F\to F$ be a diffeomorphism of a closed
surface of genus $\geq 2$. Assume that $f$ is reducible and in the JSJT
decomposition of the mapping torus of $f$ (as described in {\bf Case
3} of the proof of the Main Lemma) there are only Seifert fibered
pieces. Hence there are mutually disjoint simple closed curves $C_1,\cdots
,C_n$ on $F$ which are not null homotopic and $f(\bigcup C_i)=\bigcup
C_i$. Assume that each $C_i$ represents an element of the commutator
subgroup of $\pi_1(F)$. It is now easy to check that condition $(3)$ in
the definition of special diffeomorphism is satisfied. In fact one can
check that any component of $p^{-1}(S)$ contains a free abelian subgroup
of rank $2$ for every Seifert fibered piece $S$ of $M_f$.\endremark

\medskip
\noindent
{\bf Acknowledgement.} The author is grateful to P.B. Shalen for some
discussions regarding some parts of the proof of the Main Lemma.

\enddocument